\documentclass[]{article}
\usepackage{amssymb,enumerate,amsmath,amsthm,mathrsfs,amsfonts,graphics,graphicx,color,tikz,pinlabel}
\usepackage{float}
\usepackage[a4paper, total={6in, 8in}]{geometry}
\theoremstyle{definition}
\usepackage{indentfirst}

\title{An Algorithmic Definition of Gabai Width}
\author{Ricky Lee}
\date{}

\begin{document}

\maketitle

\newtheorem{thm}{Theorem}[section]
\newtheorem{lemma}[thm]{Lemma}
\newtheorem{prop}[thm]{Proposition}
\newtheorem{cor}[thm]{Corollary}
\newtheorem{defn}[thm]{Definition}
\newtheorem{examp}[thm]{Example}
\newtheorem{conj}[thm]{Conjecture}
\newtheorem*{question}{Question}
\newtheorem*{rmk}{Remark}
\newtheorem*{Claim}{Claim}

\begin{abstract}
	We define the Wirtinger width of a knot. Then we prove the Wirtinger width of a knot equals its Gabai width. The algorithmic nature of the Wirtinger width leads to an efficient technique for establishing upper bounds on Gabai width. As an application, we use this technique to calculate the Gabai width of approximately 50000 tabulated knots.
\end{abstract}

\section{Introduction}

Gabai width is a geometric invariant of knots that was first used by Gabai in his proof of the property R conjecture [6]. Since then, the notion of Gabai width has played central roles in many important results in 3-manifold topology such as the resolution of the knot complement problem [8], the recognition problem for $S^3$ [10], and the leveling of unknotting
tunnels [7]. The success in using Gabai width as a tool in the solution to many questions in low-dimensional topology is largely due to the deep connections between Gabai width and the topology of the knot exterior. For example, Gabai width can often be used to find incompressible surfaces ([9] Thompson, [13] Wu).

The bridge number of a knot is a closely related geometric invariant defined as the minimal number of local maxima needed to construct an embedding of the knot. A formal definition of Gabai width will follow in the subsequent section, but, roughly speaking, Gabai width depends on the number of critical points of a projection as well as their relative heights. Like most geometric invariants, both bridge number and Gabai width are notoriously difficult to calculate. However, there has been recent progress on finding algorithmically accessible definitions of bridge number. The authors in [2] defined the Wirtinger number of a link and showed that it is equal to the bridge number. The Wirtinger number is calculated using a combinatorial coloring algorithm applied to a link diagram. Using ideas inspired by the Wirtinger number, we define the Wirtinger width of a knot and show it is equal to the Gabai width of a knot.

The Wirtinger width of a knot is also calculated by coloring knot diagrams. It is algorithmically computable on a knot diagram and leads to a new and efficient combinatorial technique for establishing upper bounds on the Gabai width of a knot. In the last section of the paper, we illustrate an application of these notions by describing an algorithm used to calculate upper bounds on Gabai width, which we implemented in Python [12]. Using this implementation, we were able to calculate the Gabai width of approximately 50000 prime knots with bridge number 4, of up to 16 crossings. 

\newtheorem*{Acknowledgments}{Acknowledgments}

\begin{Acknowledgments}
	The author would like to thank Ryan Blair for introducing this topic to the author and for many helpful discussions.
\end{Acknowledgments}

\section{Preliminaries}

Let $\mathcal{K}$ denote an ambient isotopy class of knot in $\mathbb{R}^3$. Let $h:\mathbb{R}^3 \rightarrow \mathbb{R}$ defined by $h(x,y,z):=z$ be the standard height function. Let $K\subset \mathbb{R}^3$ denote a knot in the equivalence class of $\mathcal{K}$. We will always assume that the embedding of $K$ is such that $h|_K$ is a Morse function. \par 

Let $p:\mathbb{R}^3 \rightarrow \mathbb{R}^2$ defined by $p(x,y,z):=(y,z)$ be the projection map onto the yz-plane. We will always assume $K$ is embedded such that $p|_K$ is a regular projection. Then $p(K)$ is a finite four-valent graph in the yz-plane. We say that $D$ is a \emph{knot diagram} of $K$ resulting from the projection $p$ if $D$ is the graph $p(K)$ together with labels at each vertex to indicate which edges are over and which are under. By convention, these labels take the form of deleting parts of the under-arc at every crossing. Thus, we can view $D$ as a disjoint union of closed arcs in the plane. Let $\alpha_1,\ldots, \alpha_J$ denote the connected components of $D$. For each $\alpha_i$, we let $s_i$ denote the union of all edges in $p(K)$ whose interior has non-empty intersection with $\alpha_i$. We refer to each $s_i$ as a \emph{strand} and let $s(D)$ denote the set of strands of $D$. We refer to the vertices of $p(K)$ as \emph{crossings} and denote the set of vertices as $v(D)$.\par 

If $s\in s(D)$, then the two endpoints of $s$ will be referred to as the crossings \emph{incident} to $s$. If $s_p$ and $s_q$ are the under-strands of the same crossing, then we say $s_p$ and $s_q$ are \emph{adjacent}. If $\{s_1,s_2,\ldots ,s_m\}\subset s(D)$ is a subset of strands such that $s_i$ is adjacent to $s_{i+1}$ for each $i\in \{1,2,\ldots ,m-1 \}$, then we say the set $\{s_1,s_2,\ldots ,s_m\}$ is \emph{connected}. Note that there is a unique knot diagram up to planar isotopy for which there exists a strand adjacent to itself (see Figure 1). In all cases considered, we assume that adjacent strands are distinct. We say a knot diagram is \emph{trivial} if it is a diagram of the unknot. \par 

\begin{figure}[H]
	\centering	
	\includegraphics[scale=1.5]{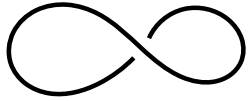}
	\caption{The unique knot diagram containing a strand adjacent to itself.}
\end{figure}

For $s\in s(D)$, we define $h(s):=max_{y\in s}h(y)$ and refer to $h(s)$ as the \emph{height of the strand} $s$. For a crossing $x\in v(D)$, we refer to $h(x)$ as the \emph{height of the crossing} $x$.
\par 

Recall that by definition, the knot diagram $D$ is a four-valent graph with labels at each vertex containing information about which edges are the under-strands. By convention, these labels take the form of deleting parts of the under-strands. However, we do not consider our labels when we calculate the height of a strand. It is therefore possible that a strand and a crossing have equal heights. In fact, if a strand is monotonic with respect to $h$, then it must have height equal to one of its incident crossings.

\begin{figure}[ht!]
	\begin{picture}(0,50)
	\put(-2,50){$s$}
	\put(-2,75){$x$}
	\end{picture}
	\centering	
	\includegraphics[scale=0.7]{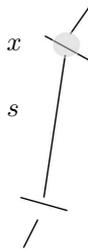}
	\caption{The strand $s$ and the incident crossing $x$ have equal heights ($h(s)=h(x)$).}
\end{figure}
\par

By \emph{critical points} of $D$ we will always be referring to images of the critical points of $h|_K$ under the projection $p$. We say that $D$ is in \emph{general position with respect to h} if all the critical points and crossings of $D$ have distinct heights with respect to $h$, $h|_K$ is Morse, and $p(K)$ is a regular projection. Observe that if the knot diagram $D$ is in general position with respect to $h$, then all the strands must have different heights.\par

We now recall the definition of Gabai width. Order the critical values of $h|_K$ by $c_1>\ldots >c_N$. Let $r_i\in (c_{i+1},c_i)$ denote arbitrarily chosen regular values of $h|_K$ for $1\leq i\leq N-1$. For any $y\in \mathbb{R}$, define $w(y):=|K\cap h^{-1}(y)|$. Define $w(K):= \sum\limits_{i=1}^{N-1}w(r_i)$. The Gabai width of $\mathcal{K}$ is defined as $\min_{K'\in \mathcal{K}}w(K')$, where the minimum is taken over all Morse embeddings of knots in the equivalence class of $\mathcal{K}$. If $K'$ is such that $w(K')=w(\mathcal{K})$, then we say $K'$ is in \emph{thin position.} \par

\section{The Coloring Rules}

In this section, we define our combinatorial method for coloring knot diagrams. Then we illustrate its connection to Gabai width. Let $D$ be a knot diagram. Let $s(D)=\{ s_1,\ldots, s_J\}$ denote the set of strands of $D$. 
\begin{defn}
	A \emph{partial coloring} is a tuple $(A,f)$ where $A$ is a subset of $s(D)$ and $f:A\rightarrow Z$ is a function with $Z\subset \mathbb{Z}$.
\end{defn}

\begin{rmk}
	Set $A_0:= \emptyset$, $Z_0:= \emptyset$, and let $f_0$ be the empty function. Then $(A_0,f_0)$ is a partial coloring. We fix $(A_0,f_0)$ to denote this vacuous partial coloring.
\end{rmk}

We define two rules for extending partial colorings. Let $(A_{t-1},f_{t-1})$ denote a partial coloring, where $t\in \mathbb{N}$ and $f:A_{t-1} \rightarrow Z_{t-1}$.

\newtheorem*{Seed Addition}{Seed Addition}

\begin{Seed Addition}
	We say the partial coloring $(A_t,f_t)$ is the result of a \emph{seed addition} to $(A_{t-1},f_{t-1})$, denoted $(A_{t-1},f_{t-1}) \rightarrow (A_t,f_t)$, if 
	
	\begin{itemize}
		\item $A_{t-1}\subset A_t$ and $A_t \setminus A_{t-1}=\{s_{i}\}$ for some strand $s_{i}\in s(D)\setminus A_{t-1}$
		\item $Z_t:=Z_{t-1}\cup \{t\}$
		\item $f_t:A_t \rightarrow Z_t$ is defined by $f_t|_{A_{t-1}}=f_{t-1}$ and $f_t(s_{i}):= t$
	\end{itemize}
	
\end{Seed Addition}

\newtheorem*{Coloring Move}{Coloring Move}

\begin{Coloring Move}
	We say $(A_t,f_t)$ is the result of a \emph{coloring move} on $(A_{t-1},f_{t-1})$, denoted $(A_{t-1},f_{t-1})\rightarrow (A_t,f_t)$, if 
	
	\begin{itemize}
		\item $A_{t-1}\subset A_t$ and $A_t \setminus A_{t-1}=\{s_q\}$ for some strand $s_q\in s(D)\setminus A_{t-1}$
		\item $s_q$ is adjacent to $s_p$ at some crossing $x\in v(D)$ and $s_p\in A_{t-1}$
		\item The over-strand $s_v$ of $x$ is an element of $A_{t-1}$
		\item $Z_t:=Z_{t-1}$
		\item $f_t:A_t \rightarrow Z_t$ is defined by $f_t|_{A_{t-1}}:= f_{t-1}$ and $f_t(s_q):=f_{t-1}(s_p)$ 
	\end{itemize}
	There are two ways we refer to a coloring move. We say that $s_q$ inherits its color from $s_p$, or that the coloring move was performed over the crossing $x$.
\end{Coloring Move}

\begin{figure}[H]
	\centering	
	\includegraphics[scale=0.5]{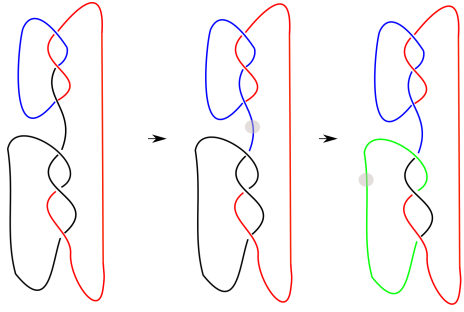}
	\caption{The transitions depict a coloring move extending the color blue, then a seed addition adding the new color green.}
\end{figure}

\begin{rmk}
	Note that we can always perform a seed addition to any uncolored strand. This allows us to use seed additions to extend the vacuous partial coloring $(A_0,f_0)$.\par 	
\end{rmk}

\begin{defn}
	If $(A_0,f_0)\rightarrow \ldots \rightarrow (A_t,f_t)$ is a sequence of coloring moves and seed additions on $D$, then we say the sequence is a \emph{partial coloring sequence}. If we have a partial coloring sequence $(A_0,f_0)\rightarrow \ldots \rightarrow (A_J,f_J)$ such that $s(D)=A_J$, then we say the sequence is a \emph{completed coloring sequence}. If $t$ is an index of a partial coloring $(A_t,f_t)$ in a specified coloring sequence, then we will refer to $t$ as a \emph{stage}.
\end{defn}

Note that we can define a completed coloring sequence for any knot diagram since we can perform a seed addition to any strand.

\begin{defn}
	If $(A_t,f_t)$ is the result of a seed addition to $(A_{t-1},f_{t-1})$ with $\{s_{i}\}=A_t\setminus A_{t-1}$, then we call $s_{i}$ a \emph{seed strand}.
\end{defn}

\begin{defn}
	Let $(A_0,f_0)\rightarrow \ldots \rightarrow (A_J,f_J)$ be a completed coloring sequence on the knot diagram $D$. Let $x\in v(D)$. Denote the over-strand of $x$ by $s_v$ and the under-strands of $x$ by $s_p$ and $s_q$. If there exists a stage $t$ such that $s_p,s_q,s_v\in A_t$ and $f_t(s_p)\neq f_t(s_q)$, then we say $x$ is a \emph{multi-colored crossing}. The smallest stage at which all previously stated conditions are satisfied will be referred to as the stage at which the crossing $x$ becomes multi-colored. 
\end{defn}

Completed coloring sequences allow us to extract geometric information from knot diagrams. To do this, we first record the order in which strands become colored, and crossings become multi-colored.

\begin{defn}
	Let $(A_0,f_0)\rightarrow \ldots \rightarrow (A_J,f_J)$ be a completed coloring sequence with multi-colored crossing set $\mathscr{C}$. Let $\mathscr{C}_t$ denote the set of crossings that become multi-colored at stage $t$. A \emph{$\Delta$-ordering} is an enumeration of the elements in $s(D)\cup \mathscr{C}$, $\Delta:=(d_i)_{i=1}^{|s(D)|+|\mathscr{C}|}$, satisfying the following conditions:
	\begin{enumerate}
		\item  For all $0\leq t< u\leq J$, all elements colored (or multi-colored) at stage $t$ are listed before any element colored (or multi-colored) at stage $u$.
		\item For each stage $0\leq t\leq J$, the element in $A_t\setminus A_{t-1}$ is listed, followed by all elements in $\mathscr{C}_t$ (if $\mathscr{C}_t\neq \emptyset$).
\end{enumerate}
	
\end{defn}

Later, we use $\Delta$-orderings to reconstruct an embedding of our knot in $\mathbb{R}^3$ from a colored knot diagram. Each seed strand will induce a single maximum and each multi-colored crossing will induce a single minimum in our reconstructed embedding. The ordering of the critical points, by decreasing height with respect to $h$, is reflected in our $\Delta$-ordering. We now show how to elevate this relationship into a calculation of Gabai width. 

\begin{defn}
 Let $(A_0,f_0) \rightarrow ... \rightarrow (A_J,f_J)$ be a completed coloring sequence. Let $\mathscr{S}\subseteq s(D)$, $\mathscr{C}\subseteq v(D)$, and $\Delta$ be the seed strands, multi-colored crossings, and $\Delta$-ordering, respectively, of our completed coloring sequence. Let $\Delta':=(d_{i_j})_{j=1}^{|\mathscr{S}|+|\mathscr{C}|}$ denote the 
 subsequence of $\Delta$ formed by restricting our $\Delta$-ordering to the set $\mathscr{S}\cup \mathscr{C}$. We define the \emph{attached sequence} $(a_i)_{i=0}^{|\Delta'|}$ to be the sequence created via the following rule:

\begin{itemize}
	\item Set $a_0:=2$
	\item If $d_{i_j}\in \Delta'$ is a seed strand and $j>1$, then set $a_j:=a_{j-1} +2$.
	\item If $d_{i_j}\in \Delta'$ is a multi-colored crossing, then set $a_j:= a_{j-1} -2$.
\end{itemize}

If the first $t$ transitions of the completed coloring contains $|S|$ total seed additions, and $|C|$ total crossings become multi-colored by stage $t$, then we say the partial coloring sequence  $(A_0,f_0)\rightarrow \ldots \rightarrow (A_t,f_t)$ induces the first $|S|+|C|$ terms of the attached sequence $(a_i)_{i=0}^{|\Delta'|}$.

\end{defn}

\begin{defn}
Define $\mathbb{W}(D):= \min\sum \limits _{i=0}^N a_i$, where the minimum is taken over all possible completed coloring sequences defined for the diagram $D$. Let $\mathbb{W}(\mathcal{K}):=\min \mathbb{W}(D)$, where the minimum is taken over all possible knot diagrams of knots in the isotopy class of $\mathcal{K}$. We define $\mathbb{W}(\mathcal{K})$ to be the \emph{Wirtinger width} of $\mathcal{K}$.
\end{defn}

\begin{rmk}
	Note that the $\Delta$-ordering resulting from a completed coloring sequence need not be unique if there exists a stage at which multiple crossings become multi-colored. However, this does not change the attached sequence or calculation of Wirtinger Width associated to each completed coloring sequence. This is because in each possible $\Delta$-ordering, the crossings which become multi-colored at the same stage must always be listed consecutively. 
\end{rmk}

 Our main theorem is the following:

\begin{thm}\label{Main Theorem}
 If $\mathcal{K}$ is an equivalence class of knot, then $\mathbb{W}(\mathcal{K})=w(\mathcal{K})$.
\end{thm}

\begin{rmk}
	The name Wirtinger Width comes from the fact, proved in [2], that the minimum number of seed additions necessary to obtain a completed coloring sequence on the knot diagram $D$, is equal to the minimum number of meridional generators needed in a Wirtinger presentation of the knot group from a diagram.
\end{rmk}

\section{Connections to the Wirtinger Number}

In this section, we prove some preliminary results that will be needed for our proof of Theorem \ref{Main Theorem}. These results are the Wirtinger width analogues of Proposition 2.2 in [2]. Let $s(D)=\{s_1,\ldots, s_J\}$ denote the strands of the knot diagram $D$.

\begin{defn}
	Let $A:=\{s_1,\ldots,s_n\}$ be a connected subset of $s(D)$, ordered by adjacency. Let $g:A \rightarrow \mathbb{Z}$. We say $g$ has a \emph{local maximum} at $s_j$ if the function $g':\{1,2,\ldots, n\} \rightarrow \mathbb{Z}$ defined by $g'(i):=g(s_i)$ has a local maximum at $j$ (that is, $g'(j)>\max\{g(j-1),g(j+1) \}$ if $1<j<n$, $g'(1)>g'(2)$, or $g(n)>g(n-1)$. 
\end{defn}

The following is an equivalent reformulation of being \emph{k-meridionally colorable}, and the main theorem, from [2].

\begin{defn}
	$D$ is \emph{k-meridionally colorable} if there exists a completed coloring sequence $(A_0,f_0)\rightarrow \ldots \rightarrow (A_J,f_J)$ containing only $k$ seed additions.
\end{defn}

\begin{thm}
	Let $\mu(\mathcal{K})$ denote the minimal $k$ such that there exists a knot diagram $D$ of a knot in the ambient isotopy class of $\mathcal{K}$ which is k-meridionally colorable. Let $\beta(\mathcal{K})$ denote the bridge number of $\mathcal{K}$. Then $\mu(\mathcal{K})=\beta(\mathcal{K})$.
\end{thm}

\begin{prop}\label{Wirtinger Prop}
	Let $(A_0,f_0)\rightarrow \ldots \rightarrow (A_J,f_J)$ be a completed coloring sequence on a knot diagram $D$. Let $\Delta:=(d_i)_{i=1}^M$ be a $\Delta$-ordering on $s(D)\cup \mathscr{C}$ induced by the completed coloring sequence on $D$. Define $h_o:\Delta \rightarrow \mathbb{Z}$ by $h_o(d_t):=-t$. Let $x\in v(D)$ be a crossing with under-strands $s_p$ and $s_q$ and over-strand $s_v$. Let $s_p$ and $s_r$ be the strands adjacent to $s_q$.
	
	\begin{enumerate}
		\item For all $u \in \{0,1,\ldots ,J\}$ and $y\in f_u(A_u)$, $f^{-1}_u(y)$ is connected.
		\item For all $y\in f_J(A_J)$, $h_o$ has a unique local maximum on $f^{-1}_J(y)$ when the set $f^{-1}_J(y)$ is ordered sequentially by adjacency.
		\item Suppose now $D$ is a non-trivial knot diagram and $f_J(s_p)=f_J(s_q)=f_J(s_r)=y$. If $k$ is such that $\{s_q\}=A_k\setminus A_{k-1}$, then we cannot have $\{s_p,s_r\}\subset A_{k-1}$.
		\item If $D$ is non-trivial and $x\notin \mathscr{C}$, then $h_o(s_v)>min\{h_o(s_p),h_o(s_q) \}$.
		\item If $D$ is non-trivial and $x\in \mathscr{C}$, then $h_o(x)<min\{h_o(s_p),h_o(s_q),h_o(s_v)\}$
	\end{enumerate}
	
\end{prop}

\begin{proof}
	(1) We induct on the stage $u$. $A_0=\emptyset$ and $f_0$ is the empty function so the claim is vacuously true for $f_0$.\par 
	
	Suppose for induction for all $u<t$ and $y\in f_u(A_u)$ we have that $f^{-1}_u(A_u)$ is connected. We will show that for all $y\in f_t(A_t)$, $f^{-1}_t(y)$ is connected. Say $\{s_{i_j}\}=A_t \setminus A_{t-1}$ and $f_t(s_{i_j})=r$. We consider two cases.\par 
	
	First, suppose $(A_t,f_t)$ is the result of a seed addition to $(A_{t-1},f_{t-1})$. By our definition of seed addition, $f^{-1}_t(r)=\{s_{i_j}\}$ and $\forall y\in f_t(A_t)\setminus \{r\}$ we have $f^{-1}_t(y)=f^{-1}_{t-1}(y)$. $f^{-1}_t(r)$ is a singleton so it is connected. $f^{-1}_{t-1}(y)$ is connected for all $y\neq r$ by our induction hypothesis.\par 
	
	Now suppose $(A_t,f_t)$ is the result of a coloring move on $(A_{t-1},f_{t-1})$. By definition of coloring move, $f^{-1}_t(r)=f^{-1}_{t-1}(r)\cup \{s_{i_j}\}$ and $s_{i_j}$ must be adjacent to a strand in $f^{-1}_{t-1}(r)$. $f^{-1}_{t-1}(r)$ is connected by our induction hypothesis so $f^{-1}_t(r)$ must also be connected. For $y\in f_t(A_t)\setminus \{r\}$, $f^{-1}_t(y)=f^{-1}_{t-1}(y)$ is also connected by our induction hypothesis. This closes induction.\\
	
	(2) We proceed by induction on the stage $u$. By definition, $A_1$ is a singleton and $f_1:A_1 \rightarrow \{1\}$. Thus $h_o$ trivially attains a unique local maximum on the set $A_1=f^{-1}(1)$.\par 
	
	Suppose for induction for all $u<t$, $h_o$ attains a unique local maximum on the set $f^{-1}_u(y)$ when ordered sequentially by adjacency, where $y\in f_u(A_u)$ is arbitrary. We claim the same holds for $f_t$. Say $\{s_{i_j}\}=A_t \setminus A_{t-1}$ and $f^{-1}_t(s_{i_j})=r$. We consider two cases.\par 
	
	First, suppose $(A_t,f_t)$ is the result of a seed addition to $(A_{t-1},f_{t-1})$. By our definition of seed addition, $f^{-1}_t(r)=\{s_{i_j}\}$ so $h_o$ trivially attains a unique local maximum on this set. For $y\in f_t(A_t)\setminus \{r\}$, $f^{-1}_t(y)=f^{-1}_{t-1}(y)$ and our claim follows from the induction hypothesis.\par  
	
	Now suppose $(A_t,f_t)$ is the result of a coloring move on $(A_{t-1},f_{t-1})$. Then there exist a strand $s_l\in A_{t-1}$ such that $f_t(s_l)=r$ and $s_l$ is adjacent to $s_{i_j}$. By definition of coloring move and $h_o$, since $\{s_{i_j}\}=A_t\setminus A_{t-1}$ and $s_l\in A_{t-1}$, then $h_o(s_{i_j})<h_o(s_l)$. Thus $s_{i_j}$ is not a local maximum in $f^{-1}_t(r)$. $f^{-1}_t(r)=f^{-1}_{t-1}(r)\cup \{s_{i_j}\}$ and $f(y)\in f_t(A_t)\setminus \{r\}$, $f^{-1}_t(y)=f^{-1}_{t-1}(y)$ so our claim follows from the induction hypothesis. This closes induction.\\
	
	(3) Suppose for contradiction $s_p,s_r\in A_{k-1}$. By assumption, $s_q\notin A_{k-1}$. By part 1 of Proposition \ref{Wirtinger Prop}, $f^{-1}_{k-1}(y)$ is connected. $f_J(s_p)=f_J(s_r)$ so we must have $\{s_p,s_r\}\subset f^{-1}_{k-1}(y)$. D is a knot diagram so this implies $s(D)\setminus \{s_q\}=f^{-1}_{k-1}(y)$. Thus $s(D)\subset f^{-1}_J(y)$ and we get that our completed coloring sequence has only a single seed strand. By Theorem 4.3, this implies $\beta(D)=1$. But the only knot with bridge number 1 is the unknot, so this contradicts the non-triviality of $D$.\\
	
	(4) Assume for contradiction $h_o(s_v)\leq min\{h_o(s_p),h_o(s_q)\}$. $D$ is non-trivial so $s_p\neq s_q$. So say without loss of generality $h_o(s_p)>h_o(s_q)$. Say $f_J(s_p)=f_J(s_q)=y$. Let $s_p$ and $s_r$ be the strands adjacent to $s_q$.\par 
	
	Let $t$ be such that $\{s_q\}=A_t\setminus A_{t-1}$. Since $h_o(s_v)\leq min\{ h_o(s_p),h_o(s_q) \}$, then $s_q$ must have inherited its color from $s_r$. $h_o(s_p)>h_o(s_q)$. Thus, $s_q,s_r\in A_{t-1}$. But $s_q\notin A_{t-1}$ and $s_p,s_q,s_r\in f^{-1}_J(y)$. This contradicts Proposition \ref{Wirtinger Prop} part 3.\\
	
	(5) This follows directly from the definition of $\Delta$-ordering. Note that a crossing is not considered to be multi-colored until its under and over-strands have all been colored.
	
\end{proof}

\section{Coloring by Height}

In this section we describe a specific procedure for coloring knot diagrams in thin position. It will be used to establish the inequality $\mathbb{W}(\mathcal{K})\leq w(\mathcal{K})$. Our goal is to obtain a coloring sequence that induces a $\Delta$-ordering which respects the ordering of the critical points of $h|_D$ by height. \par 

For the rest of this section, let $K$ be an embedding of the knot $\mathcal{K}$ in $\mathbb{R}^3$ that is in thin position with respect to $h$. Furthermore, let $K$ be such that the knot diagram $D\subset \{\text{yz-plane} \}$, resulting from the projection $p$ into the yz plane, is in general position with respect to $h$. Let $c_1>c_2>\ldots >c_N$ be the critical values of $h|_D$ ordered by decreasing height with respect to $h$. We also assume that $\mathcal{K}$ is not the ambient isotopy class of the unknot, so that $D$ is a non-trivial diagram.\par 

\begin{defn}
	Let $L$ be any knot diagram embedded in the yz-plane that is in general position with respect to $h$. Let $x\in v(L)$. Denote the under-strands of $x$ by $s_f$ and $s_r$. If $h|_{s_f}$ has a local maximum at $x$, then we say $s_f$ is \emph{ the falling strand of $x$}. If $h|_{s_r}$ has a local minimum at $x$, then we say $s_r$ is \emph{the rising strand of $x$}.  
\end{defn}

\begin{figure}[H]
	\begin{picture}(0,50)
	\put(-2,42){$s_r$}
	\put(45,0){$s_f$}
	\end{picture}
	\centering	
	\includegraphics[scale=1.0]{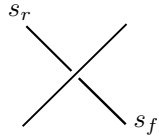}
	\caption{$s_r$ and $s_f$ denote the rising and falling strands of the pictured crossing. }
\end{figure}

\begin{defn}\label{Color by Height}
We say that we \emph{color $D$ by height} if we obtain a completed coloring sequence $(A_0,f_0)\rightarrow \ldots \rightarrow (A_J,f_J)$ by the following procedure:

\begin{enumerate}
	\item[\textbf{Step 1:}] Write $s(D)=\{s_1,\ldots, s_{|s(D)|} \}$, where $h(s_1)> \ldots > h(s_{|s(D)|})$.
	\item[\textbf{Step 2:}] Let $(A_1,f_1)$ be the result of a seed addition to $(A_0,f_0)$ such that $\{s_1\}=A_1\setminus A_0$.
	\item[\textbf{Step 3:}] Suppose we have a partial coloring sequence $(A_0,f_0)\rightarrow \ldots \rightarrow (A_{n-1},f_{n-1})$ defined, where $A_{n-1}=\{s_1,\ldots, s_{n-1} \}$. Let $x_i$ and $x_j$ be the crossings incident to $s_n$. Say $h(x_i)<h(x_j)$. We consider two cases:
	\begin{enumerate}
		\item [\textbf{Case 1:}] Suppose $h|_{s_n}$ is maximized in $int(s_n)$. Then we let $(A_n,f_n)$ be the result of a seed addition to $(A_{n-1},f_{n-1})$ such that $\{s_n\}=A_n\setminus A_{n-1}$.
		\item[\textbf{Case 2:}] Suppose $h|_{s_n}$ is maximized in $\partial s_n$ (so $s_n$ is the falling strand of $x_j$). Then we let $(A_n,f_n)$ be the result of a coloring move over $x_j$.
	\end{enumerate}	
\end{enumerate}
\end{defn}

\begin{rmk}
	It will be important to observe that when a coloring move is performed over a crossing $x$ during the color by height process, colors must extend from the rising strand of $x$ to the falling strand of $x$. Recall that in all cases considered, adjacent strands are distinct so the rising and falling strands of $x$ will always be distinct.
\end{rmk}

We first verify that knot diagrams in general position can always be colored by height.

\begin{prop}
	If $D$ is a knot diagram in general position with respect to $h$, then $D$ can be colored by height.
\end{prop}

\begin{proof}
	We verify that each step of the color by height procedure can always be performed on $D$. $D$ is in general position with respect to $h$. This means all strands have distinct heights. Thus, they can be ordered by height. By definition, we can always perform seed addition moves at any stage. What remains to be verified is that we can perform the coloring move stated in step 3 case 2.\par 
	
	Let $(A_n,f_n)$, $s_n$, $x_i$, and $x_j$ be as stated in step 3 case 2 of the color by height definition. Let $s_r$ and $s_v$ denote the over-strand and rising strand of the crossing $x_j$ respectively. Note $h(s_n)=h(x_j)$. $D$ is in general position with respect to $h$ so $h|_K$ is Morse and $p(K)$ is a regular projection. Therefore, $h(s_n)< \min\{h(s_v),h(s_r)	\}$. Thus, $\{s_v,s_r\}\subset A_{n-1}$ so we can perform the desired coloring move.
	
\end{proof}

Our goal now is to show that when we color $D$ by height, we will get $\mathbb{W}(D)\leq w(K)$. The intuition behind the following results is that since $K$ is in thin position and the resulting diagram $D$ is in general position with respect to $h$, then all strands should be in one of the following forms:

\begin{figure}[H]
	\centering	
	\includegraphics[scale=.4]{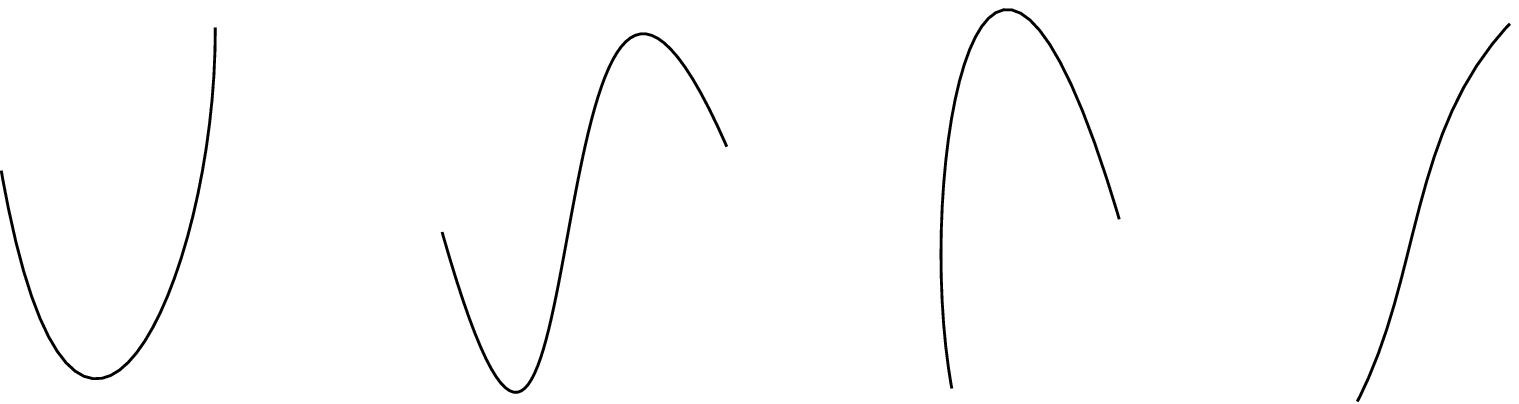}
\end{figure}

This intuition is reflected in the following Lemma.

\begin{lemma}\label{Lemma 1}
	If $s\in s(D)$ and $r\in \mathbb{R}$ is a regular value of $h|_D$, then $|s\cap h^{-1}(r)|\leq 2$.
\end{lemma}

\begin{proof}
	Suppose for contradiction we have a strand $s\in s(D)$ and a regular value $r\in \mathbb{R}$ of $h|_D$ such that $|s\cap h^{-1}(r)|\geq 3$.\par 
	
	Say $r\in (c_{j+1},c_j)$. Choose regular values $r_i\in (c_{i+1},c_i)$ for $1\leq i\leq N-1$ where $r_j=r$. Recall $K$ is in thin position, so $w(K)=w(\mathcal{K})$. To obtain our desired contradiction, we will exhibit an isotopy on $K$ to produce another embedding of $\mathcal{K}$ with strictly lower width.\par 
	
	Take three consecutive points $a,b,c$ in $s\cap h^{-1}(r)$. Let $s_{a,b}$ denote the sub-arc of $s$ in the yz plane with boundary set $\{a,b\}$. Define $s_{a,c}$ and $s_{b,c}$ similarly. Let $\alpha_{a,b}$ be the arc in $\{\text{yz plane} \} \cap h^{-1}(r)$ with boundary set $\{a,b\}$. Define $\alpha_{a,c}$ and $\alpha_{b,c}$ similarly.\par 
	
	Before describing the isotopy we must consider cases based on the order of points $\{a,b,c\}$ in $\{\text{yz plane} \}\cap h^{-1}(r)$. The ordering is by the y-coordinates of the points. Up to symmetry, there are 2 cases to consider, as depicted in the following figure.\\ 
	
	\begin{figure}[H]
		\begin{picture}(0,50)
		\put(35,30){$a$}
		\put(253,30){$a$}
		\put(75,30){$b$}
		\put(352,30){$b$}
		\put(297,30){$c$}
		\put(123,30){$c$}
		\put(50,30){$\alpha_{a,b}$}
		\put(95,30){$\alpha_{a,c}$}	
		\put(320,30){$\alpha_{b,c}$}
		\put(58,92){$s_{a,b}$}	
		\put(297,101){$s_{a,b}$}
		\put(320,-1){$s_{b,c}$}
		\put(300,42){$\alpha_{a,b}$}
		\end{picture}
		\centering	
		\includegraphics[scale=1.0]{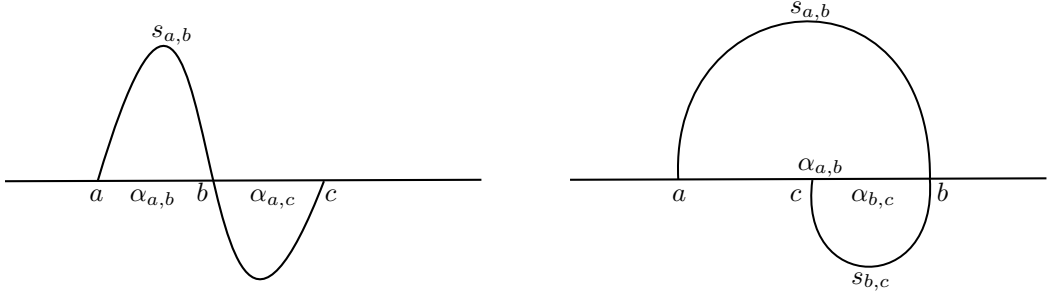}
		\caption{The set up for Case 1 and 2 are on the left and right respectively.}
		\label{StringlinkStacking}
	\end{figure}
	
	\textbf{Case 1:} Suppose $a<b<c$. Let $D_{a,b}$ be the disk cobounded by $s_{a,b}$ and $\alpha_{a,b}$. Let $D_{b,c}$ be the disk cobounded by $s_{b,c}$ and $\alpha_{b,c}$. We now define the steps of the isotopy. Let $\hat{s}_{a,c}$ be the arc component of $K\cap p^{-1}(s_{a,c})$.
	
	\begin{enumerate}
		\item[\textbf{Step 1:}] Perform an isotopy on $K$ that fixes the $y$ and $z$ coordinates of all points on $K$, arranges such that all points in $K\setminus \hat{s}_{a,c}$ have negative $x$ coordinate, and such that $\hat{s}_{a,c}=p(\hat{s}_{a,c})=s_{a,c}$. Note now $\hat{s}_{a,c}$ cobounds the two disks, $D_{a,b}$ and $D_{b,c}$ with $\alpha_{a,b}$ and $\alpha_{b,c}$ in the yz plane.
		\item[\textbf{Step 2:}] Perform an isotopy on $\hat{s}_{a,c}$ that fixes $a,b$, and $c$ and pushes $\hat{s}_{a,c}$ across $D_{a,b}$ and $D_{b,c}$ onto $\alpha_{a,c}$.
		\item[\textbf{Step 3:}] After performing the isotopy, perturb the portion of $K$ in a neighborhood of $\alpha_{a,c}$ so that $h|_K$ is Morse and has 2 fewer critical points.
	\end{enumerate}
	
	Let $s'_{a,c}$ and $K'$ denote the image of $s_{a,c}$ and $K$ respectively after the isotopy and perturbation procedure. Let $D'$ denote the diagram of $K'$ given by projection into the yz plane. Let $s'_{a,c}$ denote the image of $\hat{s}_{a,c}$ in $D'$.\\
	
	\textbf{Case 2:} Suppose $a<c<b$. Then $s_{a,c}$ cobounds a single disk with $\alpha_{a,c}$ in the yz plane and we obtain $\hat{s}'_{a,c}$, $K'$, and $D'$ from a procedure analogous to that in Step 1. The only modification is that in step 1, we push across a single disk instead of two.\\ 
	
	We now claim $w(K')<w(K)$. By construction, 
	\[	|\hat{s}'_{a,c}\cap h^{-1}(r_j)|<|\hat{s}_{a,c}\cap h^{-1}(r_j)|	\]
	Our procedure fixed the height of all points in $K$ outside of a small neighborhood of $\hat{s}_{a,c}$ and did not introduce any new critical points. Therefore,
	\[	\sum_{i=1}^{N-1}|K'\cap h^{-1}(r_i)|<\sum_{i=1}^{N-1}|K\cap h^{-1}(r_i)|=w(\mathcal{K})	\]
	The above inequality shows $w(K')<w(\mathcal{K})$. $K$ was assumed to be in thin position so we get our desired contradiction.
	
\end{proof}

\begin{prop}\label{Proposition 3}
	Let $(A_0,f_0)\rightarrow \ldots \rightarrow (A_J,f_J)$ be a completed coloring sequence obtained from coloring $D$ by height. Let $\mathscr{S}$ and $\mathscr{C}$ denote the set of seed strands and multi-colored crossings, respectively, resulting from the coloring.
	\begin{enumerate}
		\item $s\in \mathscr{S}$ if and only if $h|_s$ is maximized in $int(s)$.
		\item Let $x_i$ be a crossing with falling strand $s_q$, where $x_i$ and $x_j$ are the crossings incident to $s_q$. Then $x_i\in \mathscr{C}$ if and only if $h|_{s_q}$ is minimized in $int(s_q)$ and $h(x_i)<h(x_j)$.
	\end{enumerate}
\end{prop}

\begin{proof}
	1. By definition \ref{Color by Height}, a seed addition is performed on a strand if and only if that strand has a maximum in its interior. \\
	
	2. Let $t$ be such that $\{s_q\}=A_t \setminus A_{t-1}$.\par 
	
	($\implies$) Suppose for contradiction $h|_{s_q}$ is minimized at $x_j$. Then $h(x_i)>h(x_j)$. Note $s_q$ is the falling strand of $x_i$ and the rising strand of $x_j$. Thus, $s_q$ must be monotonic with respect to $h$. Otherwise, we could find a regular value $r$ such that $|s_q\cap h^{-1}(r)|\geq 3$, violating Lemma \ref{Lemma 1}. But this means $(A_t,f_t)$ must have been the result of a coloring move on $(A_{t-1},f_{t-1})$ over $x_i$, contradicting $x_i\in \mathscr{C}$.\par 
	
	We now show $h(x_i)<h(x_j)$. Suppose for contradiction $h(x_i)>h(x_j)$. This means if $s_q$ inherited its color via a coloring move, then the coloring move would have been performed over $x_i$ since we colored $D$ by height. But $x_i\in \mathscr{C}$ so this cannot happen. We conclude that under the current assumptions, $(A_t,f_t)$ must have been the result of a seed addition to $(A_{t-1},f_{t-1})$. By part 1, $h|_{s_q}$ must have a local maximum in $int(s_q)$. It was also assumed $h|_{s_q}$ must have a local minimum in $int(s_q)$. But this cannot happen without violating Lemma \ref{Lemma 1} since $h(x_i)>h(x_j)$ and $s_q$ is the falling strand of $x_i$. This establishes the desired contradiction.\\
	
	\begin{figure}[H]
		\begin{picture}(0,50)
		\put(83,50){$s_q$}
		\put(110,20){Level Surface}
		\end{picture}
		\centering	
		\includegraphics[scale=1.0]{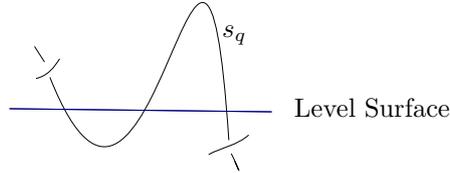}
		\caption{An example of a violation of Lemma 5.4.}
		\label{StringlinkStacking}
	\end{figure}
	
	($\impliedby$) Let $s_p$ and $s_r$ be the strands adjacent to $s_q$ at the crossings $x_i$ and $x_j$ respectively. Let $s_l$ and $s_q$ be the strands adjacent to $s_p$. Let $u$ be such that $\{s_p\}=A_u\setminus A_{u-1}$.
	
	\begin{figure}[H]
		\begin{picture}(0,50)
		\put(-6,70){$s_l$}
		\put(-6,42){$s_p$}
		\put(-5,30){$x_i$}
		\put(40,-10){$s_q$}
		\put(82,90){$s_r$}
		\put(85,67){$x_j$}
		\end{picture}
		\centering			
		\includegraphics[scale=1.0]{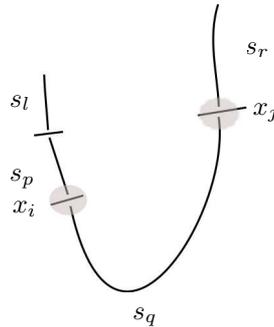}
		\caption{The setup for the current proof.}
	\end{figure}
	
	Observe that since $h(x_i)<h(x_j)$ and $s_q$ is the falling strand of $x_i$, then no coloring move could have been performed at $x_i$ when we color $D$ by height. Suppose for contradiction $x_i \notin \mathscr{C}$. We consider two cases.\par 
	
	 Recall $\{s_q\}=A_t\setminus A_{t-1}$. First, suppose $(A_t,f_t)$ was the result of a seed addition to $(A_{t-1},f_{t-1})$. $x_i \notin \mathscr{C}$ so $f_J(s_p)=f_J(s_q)$. Thus $s_p$ cannot be a seed strand. Hence, $s_p$ must have inherited its color from $s_l$ because no coloring move could have been performed over $x_i$ when we colored $D$ by height. But this means $f_J(s_l)=f_J(s_p)=f_J(s_q)$ and $\{s_l,s_q\}\subseteq A_{u-1}$ must hold. This contradicts Proposition \ref{Wirtinger Prop} part 3.\par 
	
	Now suppose $(A_t,f_t)$ was the result of a coloring move on $(A_{t-1},f_{t-1})$. No coloring move could have been performed over $x_i$ when we colored $D$ by height, so $s_q$ must have inherited its color from $s_r$. But $x_i \notin \mathscr{C}$. Therefore, $f_J(s_p)=f_J(s_q)=f_J(s_r)$. If $u<t$ (that is, if $s_p$ was colored before $s_q$), then $\{s_p,s_r\}\subset A_{t-1}$ and we have a contradiction to Proposition \ref{Wirtinger Prop} part 3.\par 
	
	Now say $t<u$ (that is, $s_q$ was colored before $s_p$). We still have $f_J(s_p)=f_J(s_q)$ so $s_p$ cannot be a seed strand under the current assumptions. Thus $s_p$ must have inherited its color from $s_l$ since no coloring move could have been performed over $x_i$ when we colored $D$ by height. This forces $f_J(s_l)=f_J(s_p)=f_J(s_q)$ and $\{s_l,s_q\}\subset A_{u-1}$, contradicting Proposition \ref{Wirtinger Prop} part 3.\par 
	
	We conclude $x_i\in \mathscr{C}$.
	
\end{proof}

\begin{cor}\label{Corollary 1}
	Recall that $K$ is in thin position and $D$ has $N$ critical points. If $\mathscr{S}$ and $\mathscr{C}$ are the sets of seed strands and multi-colored crossings resulting from a coloring of $D$ by height, then $|\mathscr{S}|+|\mathscr{C}|=N$.
\end{cor}

\begin{proof}
	Proposition \ref{Proposition 3} implies that $\mathscr{S}$ and $\mathscr{C}$ are in bijective correspondence with the set of local maxima and the set of local minima of $h|_K$ respectively. This follows because $K$ is assumed to be such that $D$ is in general position with respect to $h$. Note that since $D$ is in general position with respect to $h$, then each crossing has a unique falling strand and each strand has a lowest incident crossing.
\end{proof}

\begin{thm}\label{Coloring by Height}
	Let $(A_0,f_0)\rightarrow \ldots \rightarrow (A_J,f_J)$ be a completed coloring sequence on $D$ obtained from coloring $D$ by height. If $(a_i)_{i=0}^N$ is the attached sequence of our coloring, then $\sum\limits_{i=1}^{N}a_i \leq w(K)$.
\end{thm}

\begin{proof}
	Note that Corollary \ref{Corollary 1} verifies that the number of critical points of $K$ and the number of terms in the attached sequence $(a_i)_{i=1}^N$ resulting from our coloring of $D$ by height are equal. Let $r_n\in (c_{n+1},c_n)$ denote a regular value of $h|_D$. It suffices to show $a_n\leq w(r_n)$ for $1\leq n\leq N$. Fix one such $n$.\par  First we fix some notation. For any critical value $c_i$, let $\gamma_i$ be the unique strand at which $h^{-1}(c_i)$ fails to intersect $D$ transversely. Set $w(r_0):=0$ for notational convenience. Write 
	\[	a_n=\sum_{i=1}^na_i-a_{i-1}, \hspace{0.25cm} w(r_n)=\sum_{i=1}^n w(r_i)-w(r_{i-1})	\]
	Observe that for each $i$, we have $a_i-a_{i-1}, w(r_i)-w(r_{i-1})\in \{-2,2\}$. Thus, our goal is to show
	\begin{equation}
	\sum\limits_{i=1}^{n} a_i-a_{i-1} \leq \sum\limits_{i=1}^{n} w(r_i)-w(r_{i-1})
	\end{equation}
	Let $t$, $(A_t,f_t)$, be the stage such that $s\in A_t$ if and only if $r_n<h(s)$. We can acquire such a $t$ because our completed coloring sequence was obtained from coloring $D$ by height. Define 
	\[S_n:=\{ i | a_i-a_{i-1}=2,1\leq i\leq n\},\hspace{0.25cm}	M_n:=\{ i | w(r_i)-w(r_{i-1})=2,1\leq i\leq n\}	\] 
	
	\begin{Claim}
		$|S_n| \leq |M_n|$
	\end{Claim}
	
	\begin{proof}
	First, note that $(A_0,f_0)\rightarrow \ldots \rightarrow (A_t,f_t)$ induces at least the first $n$ terms of the sequence $(a_i)_{i=1}^N$. This is because there are $n$ critical points above $r_n$. Each strand containing a maximum above $r_n$ must have received its color via a seed addition by Proposition \ref{Proposition 3} part 1. If $\gamma_j$ contains a local minimum above $r_n$, then its lower incident crossing must become multi-colored by Proposition \ref{Proposition 3} part 2. To see why this must happen by stage $t$, recall that $D$ is in general position with respect to $h$. This means the over and under-strands of the lower incident crossing of $\gamma_j$ must have height greater than the minimum in $\gamma_j$, which in turn is greater than $r_n$.\par 
	
	Since $(A_0,f_0)\rightarrow \ldots \rightarrow (A_t,f_t)$ induces at least the first $n$ terms of the sequence $(a_i)_{i=0}^N$, then $|S_n|$ is bounded above by the number of seed additions in $(A_0,f_0)\rightarrow \ldots \rightarrow (A_t,f_t)$. But $|M_n|$ is equal to the number of such seed additions by Proposition \ref{Proposition 3} part 1 because $i\in M_n$ if and only if $c_i$ is a maximum above $r_n$. Thus $|S_n|\leq |M_n|$.
		
	\end{proof}
	
	This claim shows that the number of positive terms in $\sum\limits_{i=1}^{n} a_i-a_{i-1} $ is bounded above by the number of positive terms in $\sum\limits_{i=1}^{n} w(r_i)-w(r_{i-1})$, which verifies the inequality in equation (1).
\end{proof}

\section{Lifting a Colored Diagram}

In this section we define a procedure for retrieving an embedding of a knot from a colored knot diagram, which will just be referred to as \textit{the lifting procedure}. We want our embedding to have the following property: when the critical points are ordered by height with respect to $h$, the ordering matches the $\Delta$-ordering induced by our coloring of the diagram. It is this process that will allow us to show $ w(\mathcal{K}) \leq \mathbb{W}(\mathcal{K})$. So suppose we have a diagram $D$ of a knot in the ambient isotopy class of $\mathcal{K}$ such that $\mathbb{W}(D)=\mathbb{W}(\mathcal{K})$. Also, assume that $\mathcal{K}$ is not the ambient isotopy class of the unknot, so that $D$ is a non-trivial diagram. Let $(A_0,f_0)\rightarrow \ldots \rightarrow (A_J,f_J)$ be a completed coloring sequence on $D$ with attached sequence $(a_i)_{i=0}^N$. Let $\mathscr{S}$, $\mathscr{C}$, and $\Delta=(d_i)_{i=1}^M$ denote the set of seed strands, multi-colored crossings, and the $\Delta$-ordering on $s(D)\cup \mathscr{C}$ induced by our completed coloring sequence respectively. \par 

Recall that by definition, $D$ is a 4 valent graph with labels at each vertex to indicate which strands are over and under. In previous sections, we embedded our diagrams in the yz-plane. Now, we embed our diagram in the $z=-M-1$ plane. Furthermore, we want to view $D$ as a union of disjoint arcs in the plane. So for each $d_i\in s(D)$, let $d_i^*:=\alpha_{j_i}$, where $\alpha_{j_i}$ denotes the connected component of $D$, after deleting small sub-arcs of under-strands as dictated by the labels at each crossing, having non-empty intersection with $int(d_i)$. For $d_i\in \mathscr{C}$, let $d_i^*:=d_i$. The lifting procedure we present here is a modified version of the lifting procedure presented in [2]. Our switch in perspective on diagrams is necessary for the adaptation of the lifting in [2] to our situation.\par

There are three steps to our procedure. First we will ``lift" copies of the strands and multi-colored crossings of $D$ into planes above $D$ at heights which respect the $\Delta$-ordering of $s(D)\cup \mathscr{C}$. Then we connect our lifts at their endpoints by attaching arcs to obtain an embedded knot which orthogonally projects to $D$. Lastly, we will apply an arbitrarily small perturbation to put the knot in Morse position with critical points at the required heights.\par

We now define how to \emph{lift} a colored knot diagram. Recall that the crossings of a knot diagram are by definition just points on the plane.\par 

Embed a copy of $D$ into the plane $z=-M-1$. Define $h_o:\Delta \rightarrow \mathbb{Z}$ by $h_o(d_t):=-t$. We define the \emph{lift} of $d_t$, denoted $\hat{d_t}$, to be a copy of $d_t^*$ embedded in the plane $z=h_o(d_t)$ such that the orthogonal projection of $\hat{d_t}$ onto the $z=-M-1$ plane is $d_t^*$.\par

We now show how to connect the endpoints of the lifted strands via embedded arcs in a way that respects the crossings of $D$ upon orthogonal projection. Under-strands of multi-colored crossings will be connected via arcs that have a single minima at the lifted multi-colored crossing. This is possible because, as stated in Proposition \ref{Wirtinger Prop} part 5, multi-colored crossings were lifted to heights lower then that of its under and over-strands. Adjacent strands of the same color will be connected via monotonic arcs. This can be done in a way that respects the crossings of $D$ under orthogonal projection because of Proposition \ref{Wirtinger Prop} part 4. At such a crossing, the lifted over-strand is below at least one of the lifted under-strands.\par

Let $x$ be a crossing of $D$. Let $d_p$ and $d_q$ be the under-strands of the crossing $x$. Let $d_v$ be the over-strand of $x$. Let $\epsilon>0$ be such that the ball, denoted $B(x,\epsilon)$, in the $z=-M-1$ plane has non-empty connected intersection with the strands $d_p,d_q$ and $d_v$ and empty intersection with all other strands. Let the cylinder $B(x,\epsilon) \times \mathbb{R}$ (where $\mathbb{R}$ denotes the z-direction) also intersect $\hat{d_p}$, $\hat{d_q}$, and $\hat{d_v}$. If $x\in \mathscr{C}$, then say $x=d_i$ and observe that $B(x,\epsilon) \times \mathbb{R}$ also intersects $\hat{d_i}$. Let $B(x,\epsilon)\times \mathbb{R}$ be disjoint from all other lifts. We embed an arc connecting the lifts $\hat{d_p}$ and $\hat{d_q}$, denoted $s_{pq}$, by the following rule:

\begin{itemize}
	\item Suppose $f_J(d_p)=f_J(d_q)$. Let $s_{pq}$ be a smooth monotone arc that connects the endpoints of $\hat{d_p}$ and $\hat{d_q}$ that intersect $B(x,\epsilon)\times \mathbb{R}$. We can choose $s_{pq}$ such that it is contained in $B(x,\epsilon) \times \mathbb{R}$, disjoint from $int(\hat{d_v})$, and such that that orthogonal projection of  
	\[	(\hat{d_p} \cup s_{pq} \cup \hat{d_q} \cup \hat{d_v})\cap (B(x,\epsilon)\times \mathbb{R})	\]
	onto the $z=-M-1$ plane is $B(x,\epsilon) \cap D$, where it is understood that $s_{pq}$ projects to the deleted portions of the under-strand of $x$ in $D$.
\end{itemize}

\begin{figure}[H]
	\begin{picture}(0,50)
	\put(80,70){$\hat{d_q}$}
	\put(80,125){$\hat{d_v}$}
	\put(80,180){$\hat{d_p}$}	
	\put(50,15){$d_q$}
	\put(20,34){$d_p$}
	\put(15,10){$d_v$}
	\end{picture}
	\centering	
	\includegraphics[scale=0.50]{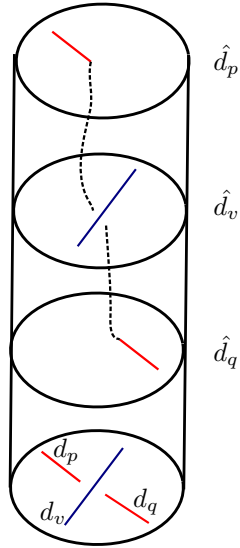}
	\caption{The construction of $s_{pq}$ (the black dashed line) at crossings that are not multi-colored. }
\end{figure}

\begin{itemize}	
	\item Suppose $f_J(d_p) \neq f_J(d_q)$. Then $d_i=x\in \mathscr{C}$. Thus $h_o(d_i)< \min \{	h_o(d_p), h_o(d_q), h_o(d_v)	\}$ by Proposition \ref{Wirtinger Prop} part 5. Let $s_{pq}$ be the union of two smooth monotone arcs connecting the endpoints of $\hat{d_p}$ and $\hat{d_q}$ in $B(x,\epsilon) \times \mathbb{R}$ to the point $\hat{d_i}$. We can choose $s_{pq}$ such that it is contained in $B(x,\epsilon) \times \mathbb{R}$, disjoint from $int(\hat{d_v})$, and such that the orthogonal projection of 
	\[(\hat{d_p} \cup s_{pq} \cup \hat{d_q} \cup \hat{d_v} \cup \hat{d_i})\cap (B(x,\epsilon)\times \mathbb{R})	\] onto the $z=-M-1$ plane is $B(x,\epsilon) \cap D$. Observe $h|_{s_{pq}}$ has a local minimum at $\hat{d_i}$ on $s_{pq}$. Again, $s_{pq}$ projects to the deleted portions of the under-strands of $x$ in $D$.
\end{itemize}  

\begin{figure}[H]
		\begin{picture}(0,50)
	\put(80,120){$\hat{d_q}$}
	\put(80,170){$\hat{d_v}$}
	\put(80,220){$\hat{d_p}$}	
	\put(50,15){$d_q$}
	\put(20,34){$d_p$}
	\put(15,10){$d_v$}
	\end{picture}
	\centering	
	\includegraphics[scale=0.5]{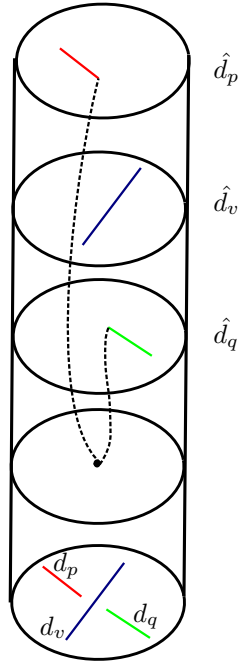}
	\caption{The construction of $s_{pq}$ (the black dashed line) at multi-colored crossings. }
\end{figure}

Now say we have lifted all strands of the diagram $D$ and performed the above step at each crossing to connect all the lifted strands. Note $\{\bigcup_t \hat{d_t}\} \cup\{ \bigcup_{p,q}s_{pq}\}\cong S^1$. However, our lifted knot is not in Morse position with respect to $h$. Our lifted strands are still parallel to the xy plane. We now define the necessary perturbation to put our knot in Morse position with respect to $h$.\par 

For each $s_{ij}$, let $y_{ij}$ denote the point in $s_{ij}$ that orthogonally projects to the corresponding crossing. Fix a strand $d_q$. Let $d_p$ and $d_r$ be strands adjacent to $d_q$. Let $[y_{pq}, y_{qr}]$ denote the sub-arc of $s_{pq} \cup \hat{d_q} \cup s_{qr}$ from $y_{pq}$ to $y_{qr}$. Recall that by Proposition \ref{Wirtinger Prop} part 2, for all $n\in f_J(A_J)$, $h_o$ has a unique local maximum on $f^{-1}_J(n)$ when the set $f^{-1}_J(n)$ is ordered sequentially by adjacency. We perturb the knot embedding via the following rule:

\begin{itemize}
	\item Suppose $d_q$ is not the unique local maximum of $h_o$ on $f^{-1}_J(f(d_q))$. Then we let the sub-arc $[y_{pq}, y_{qr}]'$ be an arbitrarily small perturbation of $[y_{pq}, y_{qr}]$ into a smooth monotonic arc, strictly increasing or decreasing as dictated by the values of $h_o(d_p)$ and $h_o(d_r)$. The perturbation is assumed to fix $y_{pq}$, $y_{qr}$, and the projection to the $z=-M-1$ plane.
	
	\item Suppose $d_q$ is the unique local maximum of $h_o$ on $f^{-1}_J(f(d_q))$. Let $m_q$ be the midpoint of $\hat{d_q}$. We let the sub-arc $[y_{pq}, y_{qr}]'$ be an arbitrarily small perturbation of $[y_{pq}, y_{qr}]$ that fixes $y_{pq}$, $m_q$, $y_{qr}$, and makes $[y_{pq}, y_{qr}]'$ strictly increase from $y_{pq}$ to $m_q$ and strictly decrease from $m_q$ to $y_{qr}$ while fixing the projection to the $z=-M-1$ plane.
\end{itemize}

\begin{rmk}
	To see why the last step, when $d_q$ is the unique local maximum of $h_o$ on $h^{-1}(h(d_q))$, is possible, it suffices to consider 2 cases. We need to verify that $h(y_{pq})<h(m_q)$ and $h(y_{qr})<h(m_q)$.\par 
	
	If $f_J(d_p)=f_J(d_q)$, then $h_o(d_p)<h_o(d_q)$ since $d_q$ is the unique local maximum of $h_o$ on $h^{-1}(h(d_q))$. Thus, $s_{pq}$ strictly decreases from the endpoint of $\hat{d_q}$ to the endpoint of $\hat{d_p}$, and $h(y_{pq})<h(m_q)$ follows.\par 
	
	If $f_J(d_p)\neq f_J(d_q)$, then $y_{pq}$ is the lift of a multi-colored crossing and the desired inequality $h(y_{pq})<h(m_q)$ follows by Proposition \ref{Wirtinger Prop} part 5. The argument for why $h(y_{qr})<h(m_q)$ holds is similar. 
\end{rmk}

We perform the stated perturbation at each lifted strand to get an embedded knot in Morse position with respect to $h$. Let $K$ denote the knot resulting from our lifting procedure.

\begin{thm}\label{Lifting}
	If $K$ is obtained from our lifting procedure on the colored knot diagram $D$, then $w(K)= \mathbb{W}(D)$. In particular, $w(\mathcal{K})\leq \mathbb{W}(\mathcal{K})$.
\end{thm}

\begin{proof}
Let $\Delta':=(d_{i_j})_{j=1}^N$ denote the subsequence of $\Delta$ formed by restricting our $\Delta$-ordering, induced by the completed coloring sequence on $D$, to the set $\mathscr{S} \cup \mathscr{C}$. Note the only critical points of $K$ are midpoints of lifted seed strands, which are maxima, and lifted multi-colored crossings, which are minima. This is because the sub-arcs of $K$ connecting these points were perturbed into monotonic arcs. Thus, $K$ has $|\mathscr{S}|+|\mathscr{C}|=N$ critical points.\par 

Let $c_1>c_2>\ldots >c_N$ be an ordering of the critical values of $h|_K$ by decreasing height. Let $1\leq j\leq N$. Since $\hat{d}_{i_j}$ is contained in the plane $z=h_o(d_{i_j})$, then $c_j=h_o(d_{i_j})$ for all $d_{i_j}\in \Delta'$. If $d_{i_j}$ is a seed strand, then $c_j=h_o(d_{i_j})$ must be a maximum. If $d_{i_j}$ is a multi-colored crossing, then $c_j=h_o(d_{i_j})$ must be a minimum. This means for any regular value $r_i\in (c_{i+1},c_i)$, $1\leq i\leq N-1$, of $h|_K$, $a_i=w(r_i)$ for any term $a_i$ in our attached sequence $(a_i)_{i=0}^N$. This shows $w(K)=\mathbb{W}(D)$, which implies $w(\mathcal{K})\leq \mathbb{W}(\mathcal{K})$.
\end{proof}

\section{Proof of Main Theorem}

We now restate and prove our main theorem.

\newtheorem*{Theorem 3.8}{Theorem 3.8}

\begin{Theorem 3.8}
If $\mathcal{K}$ is an equivalence class of knot, then $\mathbb{W}(\mathcal{K})=w(\mathcal{K})$.
\end{Theorem 3.8}

\begin{proof}
	We begin with the case where $\mathcal{K}$ is not the ambient isotopy class of the unknot. Suppose we have an embedding $K$ of a knot in the isotopy class of $\mathcal{K}$ which is in thin position with respect to $h$. Let $D$ be the knot diagram of $K$ resulting from projecting $K$ into the yz plane. By Theorem \ref{Coloring by Height}, we can color $D$ by height to see that $\mathbb{W}(\mathcal{K})\leq \mathbb{W}(D)\leq w(K)=w(\mathcal{K})$. Theorem \ref{Lifting} gives $w(\mathcal{K})\leq \mathbb{W}(\mathcal{K})$, so we get the desired equality.\par 
	
	Now suppose that $\mathcal{K}$ is the ambient isotopy class of the unknot. Then $w(\mathcal{K})=2$. We can obtain a completed coloring sequence on the standard diagram of the unknot, with no crossings, by performing a single seed addition. This shows $\mathbb{W}(\mathcal{K})\leq 2$. We now verify that $\mathbb{W}(\mathcal{K})\geq 2$. Let $U$ be a diagram resulting from an arbitrary embedding of the unknot. Let $(A_0,f_0)\rightarrow \ldots \rightarrow (A_J,f_J)$ be a completed coloring sequence on $U$ with attached sequence $(a_i)_{i=0}^N$. Let $|C_t|$ denote the number of crossings that become multi-colored at some stage $i$ with $i\leq t$. We will show $a_i\geq 0$ for $1\leq i\leq N$.\par 
	
	\begin{Claim}
		At any stage $t$, $|C_t|\leq |f_t(A_t)|$.
	\end{Claim}

	\begin{proof}
		$f^{-1}_t(y)$ is connected for all $y\in A_t$ by Proposition \ref{Wirtinger Prop} part 1. This means each set $f^{-1}_t(y)$ can contain at most two strands which are under-strands of a crossing that becomes multi-colored by stage $t$. Since every multi-colored crossing has under-strands of different colors, we get the claim, $|C_t|\leq |f_t(A_t)|$. 
	\end{proof}
	
	Now suppose for contradiction we have some $a_n < 0$ for some $1\leq n\leq N$. Without loss of generality, suppose $a_n$ is the most negative term in the attached sequence $(a_i)_{i=0}^N$ with this property. Then there exists a stage $u$ such that the first $u$ terms in our completed coloring sequence, $(A_0,f_0)\rightarrow \ldots \rightarrow (A_u,f_u)$, induces the first $n$ terms, $(a_i)_{i=0}^n$, in our attached sequence. Write 
	\[a_n= \sum_{i=1}^n a_i-a_{i-1}	\]
	Observe that 
	\[	|\{i|a_i-a_{i-1}=-2, 1\leq i \leq n \}|=|C_n|, \hspace{0.25cm} |\{i|a_i-a_{i-1}=2,1\leq i\leq n	\}|=|f_n(A_n)|	\]
	These equalities follow since the first $u$ transitions are assumed to induce the first $n$ terms of our attached sequence. Since $\sum\limits_{i=1}^n a_i-a_{i-1} =a_n<0$, then the sum must contain more negative terms then positive terms. Note $a_i-a_{i-1}\in \{-2,2\}$ by definition. This implies $|C_n|> |f_n(A_n)|$, which contradicts our claim. We conclude $a_i\geq 0$ for all $a_i$ in our attached sequence. \par 
	
	Our conclusion verifies that $\mathbb{W}(U)\geq 2$. $U$ was obtained from an arbitrarily chosen embedding of the unknot, so we get $2 \leq \mathbb{W}(\mathcal{K})$. Therefore, $\mathbb{W}(\mathcal{K})=2=w(\mathcal{K})$.
	 	 
\end{proof}

\section{Applications and Further Questions}
For an application of Theorem \ref{Main Theorem}, we focus on the tabulated knots with bridge number 4. These knots are all known to be prime. Thus, they have Gabai width equal to either 28 or 32. The algorithm we wrote and implemented takes advantage of these facts in the following way.\par 

Suppose we have a diagram of the four bridge knot $\mathcal{K}$. We want to obtain a completed coloring sequence on the diagram that starts with 3 seed additions, followed by coloring moves until we get a multi-colored crossing. Then we want to perform a single extra seed addition and finish coloring the diagram using coloring moves. If we are successful at doing this for a diagram of $\mathcal{K}$, then we know $\mathbb{W}(\mathcal{K})=w(\mathcal{K})=28$.\par 

By modifying the code in [11], which is the original algorithm for calculating Wirtinger number developed by the authors in [2], we were able to implement the above strategy to verify that approximately 50000 tabulated knots actually have Gabai width 28. The data and code for width is available at [12]. We remark that it was important to know the Gauss codes we were working on had diagrams such that the code in [11] can actually detect Wirtinger number 4. In general, this does not always happen. In [4], the authors give examples of prime, reduced, alternating diagrams $D$ of a knot $\mathcal{K}$ such that the Wirtinger number $\mu(D)$ is strictly greater then the bridge number $\beta(\mathcal{K})$.\par   

We briefly describe the verification. In [1], the authors give a method of establishing bridge number based on homomorphisms from the knot group to Coxeter groups. In the ongoing work [3], the authors use computational methods to find homomorphisms as described in [1] to verify that each of the knots tested in our code [12] have bridge number 4.

Our implementation depended heavily on the Wirtinger number of a knot diagram. In general, the search for the minimum $\mathbb{W}(D)$ over all possible diagrams $D$ is subtle. We took great advantage of the fact that the diagrams we worked on actually realized the Wirtinger number $\mu(D)$. In order to find a more robust implementation of our notions, it is important to understand how Wirtinger number and Wirtinger width interact. This leads to the following natural questions.

\begin{question}
	How can we determine whether or not a diagram $D$ realizes the minimal $\mathbb{W}(D)$ without knowing beforehand that it realizes the minimal $\mu(D)$, the Wirtinger number?
\end{question}

\begin{question}
	If the knot diagram $D$ realizes the Wirtinger number, then does $D$ also realize the Wirtinger width? 
\end{question}

One expects the answer to the second question to be no, since in [5] the authors exhibit a knot $\mathcal{K}$ such that the thin position embedding has more that $\beta(\mathcal{K})$ many maxima. However, finding a knot diagram which disproves our question seems difficult. An obvious first step is to check our knot data for a knot such that our algorithm outputs an upper bound of 32 for Gabai width, and try to show that the Gabai width of such a knot is actually 28.

\bibliography{bib}

\end{document}